\documentclass[a4paper,leqno]{article}

\usepackage{graphicx}
\usepackage{amsmath}
\usepackage{amsfonts}
\usepackage{enumerate}
\usepackage{theorem}

\usepackage{calc}

\theoremstyle{plain}

{\theorembodyfont{\rmfamily}

}

\renewcommand{\eqref}[1]{\textnormal{(\ref{#1})}}

\numberwithin{equation}{section}

\title{Reconstruction of cracks and material losses by perimeter-like penalizations and phase-field methods: numerical results}

\author{Draft}

\author{Wolfgang Ring\footnotemark[1]\hspace{1em}and\hspace{1em}Luca Rondi\footnotemark[2]}

\date{}

\begin{document}

\maketitle
\footnotetext[1]{University of Graz, Institute of Mathematics, Heinrichstr. 36 8010 Graz,
Austria. E-mail: \texttt{wolfgang.ring@uni-graz.at}}
\footnotetext[2]{Universit\`a degli Studi di Trieste, Dipartimento di Matematica e Informatica, via Valerio 12/1 34127
Trieste, Italy. E-mail: \texttt{rondi@units.it}}

\setcounter{section}{0}
\setcounter{secnumdepth}{1}

\begin{abstract}
We numerically implement the variational approach for reconstruction in the inverse crack and cavity problems developed by one of the authors. The method is based on a suitably adapted free-discontinuity problem. Its main features are the use of phase-field functions to describe the defects to be reconstructed and the use of perimeter-like 
penalizations to regularize the ill-posed problem.

The numerical implementation is based on the solution of the corresponding optimality system by a gradient method. Numerical simulations are presented to show the validity of the method.

\medskip


\medskip

\end{abstract}

\section{Introduction and setting of the method}

We consider a homogenous and isotropic conducting body, assumed to be contained in $\Omega$,
a bounded, Lipschitz domain of $\mathbb{R}^N$, $N\geq 2$. We assume that there exist $\Omega_1$, a Lipschitz domain contained in, and different from, $\Omega$, and a closed set $\gamma\subset\partial\Omega\cap\partial\Omega_1$ such that the interior of $\gamma$ is not empty and
$\gamma$ has a positive distance from $\Omega\backslash\overline{\Omega_1}$.
We assume that $\gamma$ is known and accessible to measurements.

In the body there might be present some defects, which we assume to be
perfectly insulating and outside $\Omega_1$. Namely, we model these
defects by a closed set $K_0\subset\overline{\Omega}$ such that $K_0\cap\overline{\Omega_1}$ is empty. We notice that $K_0$ represents the union of the boundaries of
these defects and that we denote with $G_{K_0}$ the connected component of $\Omega\backslash K_0$ containing $\Omega_1$, that is the region of $\Omega$ reachable from $\gamma$ without crossing $K_0$.

The defects may have different geometrical properties. For instance, we may have, even at the same time,
\emph{cracks} (either interior or surface-breaking), or \emph{material losses} (either interior, that is cavities, or at the boundary).
We recall that a defect $K_0$ is a \emph{material loss} if $G_{K_0}$ coincides with the interior of its closure.

Let us consider the following experiment.
If a current density $f_0$ is applied on $\gamma$,
then
the electrostatic potential in $\Omega$, $u_0=u(f_0,K_0)$,
is the solution to the following (normalized) Neumann boundary value problem
\begin{equation}\label{Dirpbm}
\left\{\begin{array}{ll}
\Delta u=0&\text{in }\Omega\backslash K_0\\
\nabla u\cdot \nu=f_0&\text{on }\gamma\\
\nabla u\cdot \nu=0&
\text{on }\partial(\Omega\backslash K_0)\backslash\gamma\\
\int_{\gamma}u=0.&
\end{array}
\right.
\end{equation}
The current density is modeled by a function 
$f_0\in L^s(\gamma)$, for some constant $s>N-1$, such that
$\int_{\gamma}f_0=0$.
The electrostatic potential $u_0$ may then be measured on $\gamma$. We call such a measurement $g_0=u_0|_{\gamma}$ and we observe that
$g_0\in L^2(\gamma)$ and $\int_{\gamma}g_0=0$.
In this way we obtain an electrostatic boundary measurement of voltage, $g_0$, and current, $f_0$, type on $\gamma$. In mathematical words, we measure the Cauchy data $(g_0,f_0)$ of the harmonic function $u_0$ on $\gamma$.
Clearly, prescribed the current $f_0$, the voltage $g_0$ depends on $K_0$.
If $K_0$ is unknown, then the measured voltage $g_0$ may provide information about the unknown defect.
In fact,
the aim of the inverse problem is to reconstruct an unknown defect $K_0$ by prescribing one or more current densities $f_0$ and measuring the corresponding value of the potentials on $\gamma$. Such a problem arises, for instance, in non-destructive evaluation, for the determination of flaws like cracks or cavities in conducting bodies by
non-invasive methods. We refer to this problem as the \emph{inverse crack problem}, in the general case. Instead, when we a priori know that the defect is a material loss, we denote it as the \emph{inverse cavity problem}.
For results on the inverse crack problem and related problems, we refer to the review article \cite{Bry e Vog04}. Here we simply wish to note that a single measurement (that is performing the experiment previously described only once) is enough to determine uniquely a material loss. In the general crack case, instead, one measurement may not be enough, however two suitably chosen measurements (corresponding to two suitable prescribed current densities) are enough for unique identification of any kind of defects
at least in the planar case.

Let us remark here that if the unknown defect is a priori assumed to be interior (that is $K_0\subset\Omega)$ and if the whole boundary of $\Omega$ is accessible, then
we may simply take $\gamma=\partial\Omega$.

Our approach to this inverse problem is the following.
We observe that $u_0$ is smooth outside $K_0$, whereas it may, and generally does, jump across $K_0$. Therefore,
starting from the Cauchy data, we wish to reconstruct the function $u_0$ in $\Omega$,
and in particular its discontinuity set $J(u_0)$. We notice that this is not a classical Cauchy problem for $u_0$, since $u_0$ is harmonic in $\Omega\backslash K_0$ with $K_0$ unknown! Rather, it looks more like a free-discontinuity problem for $u_0$, since
its discontinuity set $J(u_0)$ is unknown and it is actually the aim of our reconstruction.
If we are able to reconstruct $u_0$ and $J(u_0)$, then we obtain valuable information on $K_0$, given the fact the $J(u_0)\subset K_0$. Actually, for the inverse cavity problem,
$J(u_0)$ determines the whole $\partial G_{K_0}$. On the contrary, in the inverse crack problem, it may happen that a crack is not visible for a particular measurement, that is 
$J(u_0)$ does not detect the whole $\partial G_{K_0}$. In this case, we may change the
prescribed current density, reconstruct  again the electrostatic potential from its values on $\gamma$, and recover another portion of
$\partial G_{K_0}$. The uniqueness results tell us how many times and with which kind
of prescribed current densities we need to repeat this procedure to fully reconstruct the unknown defect.

The main difficulties in the reconstruction of $u_0$ from its Cauchy data are the following. First of all, the problem is severely ill-posed, as Cauchy problems for elliptic equations are. Second, since the potential $u_0$ to be reconstructed is a discontinuous function whose discontinuities are unknown (actually they are the aim of our reconstruction), the problem is not even linear. Thus all the main difficulties of the original inverse problem are still present in the reconstruction of $u_0$.

The way to tackle ill-posedness is crucial. In fact, since the boundary data are measured,
the data which are really available are not the exact Cauchy data $(g_0,f_0)$ but some noisy perturbation of them. Namely, the available data we assume to know are
$(g_{\varepsilon},f_{\varepsilon})$. Here $f_{\varepsilon}$ belongs to $L^s(\partial\Omega)$ and satisfies
$\mathrm{supp}(f_{\varepsilon})\subset\gamma$ and $\int_{\partial\Omega}f_{\varepsilon}=0$, whereas
$g_{\varepsilon}$ belongs to $L^2(\gamma)$ and satisfies $\int_{\gamma}g_{\varepsilon}=0$.
We assume that
\begin{equation}\label{noiselevel}
\|f_0-f_{\varepsilon}\|_{L^s(\gamma)}\leq \varepsilon \quad\text{and}\quad
\|g_0-g_{\varepsilon}\|_{L^2(\gamma)}\leq \varepsilon.
\end{equation}
where $\varepsilon$, $0<\varepsilon\leq 1/2$, denotes the noise level.

As mentioned, rather than a classical Cauchy problem, we consider such a  reconstruction 
as a free-discontinuity problem for the unknown potential $u_0$.
We follow the variational approach developed in \cite{Ron08,Ron09} for cracks and material losses, respectively. Such a method is based on the following two features. The first one is the choice of the regularization. In order to regularize the problem a perimeter-like penalization is used. Namely, we penalize the $(N-1)$-dimensional measure of the
unknown defect $K_0$ (actually of the discontinuity set of the unknown potential).
Second, we model discontinuity sets through phase-field functions, thus obtaining a formulation in which a discontinuous function $u$ and its discontinuity set $J(u)$ are replaced, respectively, by a smooth function $u$ and by a smooth phase-field function 
$v$. Such a formulation is amenable to numerical implementation.

In particular, for the crack case, in \cite{Ron08} it has been used a regularization
based on the so-called Mumford-Shah functional, \cite{Mum e Sha89}, and its approximation, in the sense of  $\Gamma$-convergence, with phase-field functionals
due to Ambrosio and Tortorelli, \cite{Amb e Tor1,Amb e Tor2}. For material losses, \cite{Ron09}, it has been used a more classical perimeter penalization and its approximation, again in the sense of  $\Gamma$-convergence, with phase-field functionals
due to Modica and Mortola, \cite{Mod e Mor}. For further details on free-discontinuity problems and their approximations we refer for instance to \cite{Amb e Fus e Pal,Bra}.

In this paper we develop the numerics of the approach in \cite{Ron08,Ron09}. The convergence analysis done in these papers provides a justification of the numerical method, in particular for the material loss case,
see \cite[Theorem~4.2]{Ron09}. For what concerns the crack case, we do not have a precise convergence result for the method implemented here.
However this is quite simpler from a numerical point of view
than the one developed in \cite{Ron08} and for which we have convergence results. Moreover we believe that this simplification
might still lead to good reconstructions, see also the discussion in Section~5 of \cite{Ron09}.

We shall use the following notation.
We fix a constant $q>2$, depending on some regularity properties of the defect $K_0$ to be reconstructed. Namely, we assume that there exist a constant $q>2$  and a constant $C$, independent of $f_0$, such that
$$\|\nabla u_0\|_{L^q(\Omega)}\leq C\|f_0\|_{L^s(\gamma)},$$
where $u_0=u(f_0,K_0)$ solves \eqref{Dirpbm}.
We also set
$$0<q_1=(q-2)/(2q)<1/2.$$

The function $\psi:\mathbb{R}\to\mathbb{R}$ is continuous and non-decreasing and such that $\psi(0)=0$, $\psi(t)>0$ if $t>0$, and $\psi(1)=1$.
Then, for $0<\varepsilon\leq 1/2$, we define
$$\psi_{\varepsilon}=(1-\varepsilon^2)\psi+\varepsilon^2.$$

We introduce a single-well potential $V$ centered at $1$, that is a non-negative continuous function such that $V(t)=0$ if and only if $t=1$.

We shall also need a double-well potential $W$ centered at $0$ and $1$,
that is a non-negative continuous function such that $W(t)=0$ if and only if $t\in\{0,1\}$.

We assume that the functions $W$, $V$ and
$\psi$ are $C^1$ and are bounded
all over $\mathbb{R}$ and also their derivatives are bounded and uniformly continuous all over $\mathbb{R}$.
We shall also assume that $W'$, $V'$ and $\psi'$ are H\"older continuous for some exponent $\alpha$, $0<\alpha\leq 1$,
all over $\mathbb{R}$.
Furthermore, the following assumption will be made. About $\psi$
we require that for any $t\leq 0$ we have $\psi(t)=\psi(0)=0$ and
$\psi(t)=\psi(1)=1$ for any $t\geq 1$. In particular, we have that $\psi'(0)=\psi'(1)=0$.
We also require that for any $t\leq 0$
we have $V(t)\geq V(0)$. Obviously, we have that, for any $t\leq 0$, $W(t)\geq W(0)$, and,
for any $t\geq 1$, $V(t)\geq V(1)$ and $W(t)\geq W(1)$.

For example, the following choices may be made. For any $t\in[0,1]$
$$\psi(t)=-2t^3+3t^2,\quad
W(t)=9t^2(t-1)^2,\quad
V(t)=(t-1)^2/4,$$
with straightforward extensions beyond $[0,1]$.

We define the space $W(\Omega)=\{\tilde{v}\in W^{1,2}(\Omega):\ \tilde{v}=0\text{ a.e. in }\Omega_1\}$. To any $\tilde{v}\in W(\Omega)$ we associate the function $v=1-\tilde{v}$. We remark that
$v\in W^{1,2}(\Omega)$ and $v=1$ almost everywhere in $\Omega_1$.
We finally fix positive tuning parameters $a$, $b$ and $c$, and
a noise level $\varepsilon$, $0<\varepsilon\leq 1/2$.
All these constants and the notation will be kept fixed throughout the paper.

Basically, the method is the following. Beginning from the crack case, we wish to minimize, with respect to the phase-field variable $\tilde{v}\in W(\Omega)$, with the constraint $0\leq \tilde{v}\leq 1$, the functional
$\mathcal{F}_{\varepsilon}:W(\Omega)\to \mathbb{R}$, which is defined 
as follows. For any $\tilde{v}\in W(\Omega)$, recalling that $v=1-\tilde{v}$, we set
\begin{equation}
\mathcal{F}_{\varepsilon}(\tilde{v})=
\frac{a}{\varepsilon^{q_1}}\int_{\gamma}|\tilde{u}_{\varepsilon}-g_{\varepsilon}|^2+
\displaystyle{b\int_{\Omega}\psi_{\varepsilon}(v)|\nabla \tilde{u}_{\varepsilon}|^2+\frac{c^2}{\varepsilon}
\int_{\Omega}V(v)+\varepsilon \int_{\Omega}|\nabla v|^2}.
\end{equation}
Here
$\tilde{u}_{\varepsilon}=\tilde{u}_{\varepsilon}(\tilde{v})$ solves
\begin{equation}
\left\{\begin{array}{ll}
\mathrm{div}(\psi_{\varepsilon}(v)\nabla \tilde{u}_{\varepsilon})=0 &\text{in }\Omega\\
\psi_{\varepsilon}(v)\nabla \tilde{u}_{\varepsilon}\cdot\nu=f_{\varepsilon} &\text{on }\partial\Omega\\
\int_{\gamma}\tilde{u}_{\varepsilon}=0.
\end{array}\right.
\end{equation}

We notice that the first term is the fidelity term with respect to the measured boundary datum, the other three terms are the Ambrosio-Tortorelli functional. The link with the prescribed boundary datum and  with the presence of cracks is through $\tilde{u}_{\varepsilon}$, the solution of the weighted elliptic equation. We observe that the single-well potential $V$ forces the phase-field function $v=1-\tilde{v}$ to be equal to $1$ except in a small region, which is where the crack should be located.
The tuning parameters $a$, $b$ and $c$ allow to put more emphasis on one or the other of the features of the functional. Namely, $a$ controls the match with the Dirichlet datum, $b$ the smoothness of the reconstructed potential away from its discontinuities and $c$ the penalization on the $(N-1)$-dimensional measure of the discontinuities. Therefore $c$ may be seen as a regularization parameter.

For the material loss case, we simply replace the
single-well potential $V$ with the double-well potential $W$.
Namely,
we define $\mathcal{G}_{\varepsilon}:W(\Omega)\to \mathbb{R}$ in an analogous way by simply replacing $V$ with $W$, that
is,
for any $\tilde{v}\in W(\Omega)$, we set
\begin{equation}
\mathcal{G}_{\varepsilon}(\tilde{v})=
\frac{a}{\varepsilon^{q_1}}\int_{\gamma}|\tilde{u}_{\varepsilon}-g_{\varepsilon}|^2+
\displaystyle{b\int_{\Omega}\psi_{\varepsilon}(v)|\nabla \tilde{u}_{\varepsilon}|^2+\frac{c^2}{\varepsilon}
\int_{\Omega}W(v)+\varepsilon \int_{\Omega}|\nabla v|^2}.
\end{equation}
We then minimize, with respect to the phase-field variable $\tilde{v}\in W(\Omega)$, with the constraint $0\leq \tilde{v}\leq 1$, the functional $\mathcal{G}_{\varepsilon}$.

We notice that in this case the last two terms are the Modica-Mortola functional, which penalizes the perimeter of $G_{K_0}$ in $\Omega$. Here the double-well potential $W$ forces the phase-field function $v$ to be either $0$ (inside the material loss) or $1$ (outside the material loss), with a quick transition between these two regions.

Summarizing, 
we shall minimize the functional $\mathcal{F}_{\varepsilon}$, when we aim to reconstruct defects such as cracks, 
and the functional $\mathcal{G}_{\varepsilon}$, when we aim to reconstruct material losses. Namely,
we wish to solve numerically the following minimization problems (depending on the properties of the unknown defect $K_0$)

\begin{enumerate}[(i)]

\item $\min\mathcal{F}_{\varepsilon}$ on $W(\Omega)$, with the constraint $0\leq \tilde{v}\leq 1$, if
$K_0$ contains portions of cracks.

\item $\min\mathcal{G}_{\varepsilon}$ on $W(\Omega)$, with the constraint $0\leq \tilde{v}\leq 1$, if
$K_0$ is a material loss defect. 

\end{enumerate}

Let us notice that, by the direct method, these minimum problems admit a solution. 

About the numerical method, in order to find the minimizers, we formulate the
corresponding optimality system and we use a gradient method, see Section~\ref{optsec} for details.
In Section \ref{numsec} numerical simulations are presented for both the single- and double-well approximations. Numerical experiments are performed for various types of defects with noise-free and noisy data-sets. 

\subsection{Acknowledgments}
The second author is partially supported by GNAMPA under 2008 and 2009 projects and by the Italian Ministry of University and Research under PRIN 2008 project. Part of this work was done during a visit of the first author to the University of Trieste, supported by GNAMPA under 2008 project, and during a visit of the second author at the University of Graz, supported by the Special Research Center ``Mathematical Optimization and Applications in Biomedical Sciences''. The authors wish also to thank Alfio Borz\`\i{} for useful discussions.

\section{Optimality system and the gradient method}\label{optsec}

We now look towards the numerical implementation of the method. 
We begin by recalling the differentiability properties of the functionals
$\mathcal{F}_{\varepsilon}$ and $\mathcal{G}_{\varepsilon}$, which have been investigated in 
\cite[Section~6]{Ron09}.

We define the following spaces. For any $p$, $2\leq p\leq +\infty$,
let us call $L_p(\Omega)=\{\tilde{v}\in L^{p}(\Omega):\tilde{v}=0\text{ a.e. in }\Omega_1\}$
and $W_p(\Omega)= W^{1,2}(\Omega)\cap L_{p}(\Omega)$, with norm $\|\tilde{v}\|_{L_p(\Omega)}=\|\tilde{v}\|_{L^{p}(\Omega)}$ and
$\|\tilde{v}\|_{W_p(\Omega)}=\|\tilde{v}\|_{L^{p}(\Omega)}+\|\nabla \tilde{v}\|_{L^{2}(\Omega)}$.
To any $\tilde{v}\in L^{2}(\Omega)$ we as usual associate the function $v=1-\tilde{v}$.
If $\tilde{v}$ belongs either to $L_p(\Omega)$ or to $W_p(\Omega)$, then $v\in L^p(\Omega)$,
$v=1$ almost everywhere in $\Omega_1$, and, provided $0\leq \tilde{v}\leq 1$
almost everywhere in $\Omega$,
we also have $0\leq v\leq 1$ almost everywhere in $\Omega$.
We observe that $W_2(\Omega)=W(\Omega)$ as previously defined.

For any $q$, $q\geq 2$, we define
$${W^{1,q}_{\gamma}(\Omega)}=\left\{u\in W^{1,q}(\Omega):\ \int_{\gamma}u=0\right\}.$$
We observe that, by a generalized Poincar\'e inequality, on 
${W^{1,q}_{\gamma}(\Omega)}$ the usual $W^{1,q}(\Omega)$ norm and the norm
$\|u\|_{{W^{1,q}_{\gamma}(\Omega)}}=\|\nabla u\|_{L^q(\Omega)}$ are equivalent. Therefore, we shall 
set this second one as the natural norm of ${W^{1,q}_{\gamma}(\Omega)}$.

We define $\mathcal{H}_{\varepsilon}:L_2(\Omega)\to W^{1,2}_{\gamma}(\Omega)$ as follows
$$\mathcal{H}_{\varepsilon}(\tilde{v})=\tilde{u}_{\varepsilon}(\tilde{v})\quad\text{for any }\tilde{v}\in L_2(\Omega).$$

There exist constants $p(\varepsilon)\geq 2$ and $q(\varepsilon)>2$, depending on $\varepsilon$ and $\alpha$, such that all the following results hold.

First, $\mathcal{H}_{\varepsilon}:L_2(\Omega)\to W^{1,q(\varepsilon)}_{\gamma}(\Omega)$, with bounded image in $W^{1,q(\varepsilon)}_{\gamma}(\Omega)$, and,
for any $\tilde{v}_0\in L_2(\Omega)$, such an operator
$\mathcal{H}_{\varepsilon}$ is differentiable in $\tilde{v}_0$ with respect to
the $L^p(\Omega)$, with $p\geq p(\varepsilon)$, and $W^{1,q(\varepsilon)}_{\gamma}(\Omega)$ norms.
Let $D\mathcal{H}_{\varepsilon}(\tilde{v}_0):L_p(\Omega)\to
W^{1,q(\varepsilon)}_{\gamma}(\Omega)$ be the differential in $\tilde{v}_0$. Then
for any $\tilde{v}$ in $L_p(\Omega)$ we have
$$D\mathcal{H}_{\varepsilon}(\tilde{v}_0)[\tilde{v}]=U_{\varepsilon}(\tilde{v}_0,\tilde{v})$$
where $U_{\varepsilon}=U_{\varepsilon}(\tilde{v}_0,\tilde{v})\in {W^{1,2}_{\gamma}(\Omega)}$ solves
the following problem
\begin{equation}\label{der}
\left\{\begin{array}{ll}
\mathrm{div}(\psi_{\varepsilon}(v_0)\nabla U_{\varepsilon})=
\mathrm{div}(\psi'_{\varepsilon}(v_0)\tilde{v}\nabla (\mathcal{H}_{\varepsilon}(\tilde{v}_0)))
 &\text{in }\Omega\\
\psi_{\varepsilon}(v_0)\nabla U_{\varepsilon}\cdot\nu=0 &\text{on }\partial\Omega.
\end{array}\right.
\end{equation}
Obviously, $v_0=1-\tilde{v}_0$. We recall that for any vector valued function $G\in L^2(\Omega,\mathbb{R}^N)$,
$\mathrm{div}(G)$ defines a functional on $W^{1,2}(\Omega)$ in the following way
$$\mathrm{div}(G)[\phi]=-\int_{\Omega}G\cdot \nabla \phi\quad\text{for any }\phi\in W^{1,2}(\Omega).$$
Therefore,
the weak formulation of \eqref{der} is looking for a function
$U_{\varepsilon}\in{W^{1,2}_{\gamma}(\Omega)}$ such that
$$\int_{\Omega}\psi_{\varepsilon}(v_0)\nabla U_{\varepsilon}\cdot\nabla \varphi
=
\int_{\Omega}\psi'_{\varepsilon}(v_0)\tilde{v}\nabla(\mathcal{H}_{\varepsilon}(\tilde{v}_0))
\cdot \nabla \varphi\quad\text{for any }
\varphi\in {W^{1,2}(\Omega)}.$$

Here, and analogously in the sequel, the differentiability has to be understood in the following sense. For any
$\tilde{v}$ in $L_p(\Omega)$
$$\mathcal{H}_{\varepsilon}(\tilde{v}_0+\tilde{v})=\mathcal{H}_{\varepsilon}(\tilde{v}_0)+
D\mathcal{H}_{\varepsilon}(\tilde{v}_0)[\tilde{v}]+R(\tilde{v})$$
where
$$\lim_{\|\tilde{v}\|_{L^p(\Omega)}\to 0}\frac{\|R(\tilde{v})\|_{W^{1,q(\varepsilon)}_{\gamma}(\Omega)}}{\|\tilde{v}\|_{L^p(\Omega)}}=0.$$

We conclude that, for any $\tilde{v}_0\in W(\Omega)$, $\mathcal{F}_{\varepsilon}$ and 
$\mathcal{G}_{\varepsilon}$ are differentiable in $\tilde{v}_0$ with respect to the $W_p(\Omega)$ norm,
with $p\geq p(\varepsilon)$. Let $D\mathcal{F}_{\varepsilon}(\tilde{v}_0),\ D\mathcal{G}_{\varepsilon}(\tilde{v}_0):W_p(\Omega)\to\mathbb{R}$
be the differentials in $\tilde{v}_0$ of $\mathcal{F}_{\varepsilon}$ and $\mathcal{G}_{\varepsilon}$, respectively.
Then, for any $\tilde{v}\in W_p(\Omega)$ we have
\begin{multline}\label{differentialf}
D\mathcal{F}_{\varepsilon}(\tilde{v}_0)[\tilde{v}]=
\frac{2a}{\varepsilon^{q_1}}\int_{\gamma}(\mathcal{H}_{\varepsilon}(\tilde{v}_0)-g_{\varepsilon})U_{\varepsilon}(\tilde{v}_0,\tilde{v})+\\
b \int_{\Omega}\left(2\psi_{\varepsilon}(v_0)\nabla \mathcal{H}_{\varepsilon}(\tilde{v}_0)\cdot\nabla U_{\varepsilon}(\tilde{v}_0,\tilde{v})-
\psi'_{\varepsilon}(v_0)|\nabla \mathcal{H}_{\varepsilon}(\tilde{v}_0)|^{2}\tilde{v}\right)+\\
\frac{c^2}{\varepsilon}\int_{\Omega}(-V'(v_0)\tilde{v})+2\varepsilon
\int_{\Omega}\nabla\tilde{v}_0\cdot\nabla \tilde{v}
\end{multline}
and
\begin{multline}\label{differentialg}
D\mathcal{G}_{\varepsilon}(\tilde{v}_0)[\tilde{v}]=
\frac{2a}{\varepsilon^{q_1}}\int_{\gamma}(\mathcal{H}_{\varepsilon}(\tilde{v}_0)-g_{\varepsilon})U_{\varepsilon}(\tilde{v}_0,\tilde{v})+\\
b \int_{\Omega}\left(2\psi_{\varepsilon}(v_0)\nabla \mathcal{H}_{\varepsilon}(\tilde{v}_0)\cdot\nabla U_{\varepsilon}(\tilde{v}_0,\tilde{v})-
\psi'_{\varepsilon}(v_0)|\nabla \mathcal{H}_{\varepsilon}(\tilde{v}_0)|^{2}\tilde{v}\right)+\\
\frac{c^2}{\varepsilon}\int_{\Omega}(-W'(v_0)\tilde{v})+2\varepsilon
\int_{\Omega}\nabla\tilde{v}_0\cdot\nabla \tilde{v}.
\end{multline}

An important remark is the following. If $N=2$, then we may actually choose $p(\varepsilon)=2$,
and we observe that $W_2(\Omega)$
is a Hilbert space, with the scalar product $\int_{\Omega}\nabla \tilde{v}_1\cdot\nabla \tilde{v}_2$ for any
$\tilde{v}_1$, $\tilde{v}_2\in W_2(\Omega)$.
If $N>2$, then it might happen that $p(\varepsilon)>2$ and that $W_{p(\varepsilon)}(\Omega)$ has not a Hilbert space structure anymore. However, since $p(\varepsilon)$ is finite, $W_{p(\varepsilon)}(\Omega)$ is still
a strictly convex real reflexive Banach space.

In the sequel we shall fix $p=p(\varepsilon)$, (with $p(\varepsilon)=2$ if $N=2$) and we
call $\mathcal{MF}_{\varepsilon}$ the following functional, which is defined on
$W^{1,2}_{\gamma}(\Omega)\times W_{p(\varepsilon)}(\Omega)$,
\begin{multline}
\mathcal{MF}_{\varepsilon}(u,\tilde{v})=
\frac{a}{\varepsilon^{q_1}}\int_{\gamma}|u-g_{\varepsilon}|^2+
\displaystyle{\int_{\Omega}\big(b\psi_{\varepsilon}(v)|\nabla u|^2+\frac{c^2}{\varepsilon}V(v)+
\varepsilon |\nabla v|^2\big)},\\
\text{for any }(u,\tilde{v})\in W^{1,2}_{\gamma}(\Omega)\times W_{p(\varepsilon)}(\Omega).
\end{multline}
Such a functional is finite for any $(u,\tilde{v})\in W^{1,2}_{\gamma}(\Omega)\times W_{p(\varepsilon)}(\Omega)$.
By similar reasonings, for any $(u_0,\tilde{v}_0)\in W^{1,q(\varepsilon)}_{\gamma}(\Omega)\times W_{p(\varepsilon)}(\Omega)$,
we have that $\mathcal{MF}_{\varepsilon}$
is differentiable in $(u_0,\tilde{v}_0)$
and
for any $(u,\tilde{v})\in W^{1,2}_{\gamma}(\Omega)\times W_{p(\varepsilon)}(\Omega)$ we have
\begin{multline}\label{differentialmf}
D\mathcal{MF}_{\varepsilon}(u_0,\tilde{v}_0)[(u,\tilde{v})]=
\frac{2a}{\varepsilon^{q_1}}\int_{\gamma}(u_0-g_{\varepsilon})u+\\
b \int_{\Omega}\left(2\psi_{\varepsilon}(v_0)\nabla u_0\cdot\nabla u-
\psi'_{\varepsilon}(v_0)|\nabla u_0|^{2}\tilde{v}\right)+
\frac{c^2}{\varepsilon}\int_{\Omega}(-V'(v_0)\tilde{v})+2\varepsilon
\int_{\Omega}\nabla\tilde{v}_0\cdot\nabla \tilde{v}.
\end{multline}

We observe that $\mathcal{F}_{\varepsilon}(\tilde{v})=\mathcal{MF}_{\varepsilon}(\mathcal{H}_{\varepsilon}(\tilde{v}),\tilde{v})$.
Analogously, we define $\mathcal{MG}_{\varepsilon}$ simply by replacing $V$ with $W$. Analogous properties of differentiability
hold for $\mathcal{MG}_{\varepsilon}$ as well.

Let us finally define
$\mathcal{LF}_{\varepsilon}:
W^{1,2}_{\gamma}(\Omega)\times W_{p(\varepsilon)}(\Omega)\times W^{1,2}(\Omega)\to\mathbb{R}$ 
such that for any $(u,\tilde{v},\phi)\in W^{1,2}_{\gamma}(\Omega)\times W_{p(\varepsilon)}(\Omega)\times W^{1,2}(\Omega)$
we have
\begin{equation}\label{lagrangian}
\mathcal{LF}_{\varepsilon}(u,\tilde{v},\phi)=\mathcal{MF}_{\varepsilon}(u,\tilde{v})+\int_{\Omega}\psi_{\varepsilon}(v)\nabla u\cdot \nabla\phi-\int_{\partial\Omega}f_{\varepsilon}\phi.
\end{equation}
In an analogous way we define $\mathcal{LG}_{\varepsilon}$ replacing $\mathcal{MF}_{\varepsilon}$ with
$\mathcal{MG}_{\varepsilon}$.

We observe that $\mathcal{LF}_{\varepsilon}$ (and $\mathcal{LG}_{\varepsilon}$ as well) is differentiable in any $(u_0,\tilde{v}_0,\phi_0)
\in W^{1,q(\varepsilon)}_{\gamma}(\Omega)\times W_{p(\varepsilon)}(\Omega)\times W^{1,2}(\Omega)$. For any $(u,\tilde{v},\phi)\in
W^{1,2}_{\gamma}(\Omega)\times W_{p(\varepsilon)}(\Omega)\times W^{1,2}(\Omega)$ we have
\begin{multline}\label{duweak}
\frac{\partial \mathcal{LF}_{\varepsilon}}{\partial u}(u_0,\tilde{v}_0,\phi_0)[u]=\\
\frac{2a}{\varepsilon^{q_1}}
\int_{\gamma}(u_0-g_{\varepsilon})u+
2b\int_{\Omega}\psi_{\varepsilon}(v_0)\nabla u_0\cdot\nabla u+ \int_{\Omega}\psi_{\varepsilon}(v_0)\nabla \phi_0\cdot \nabla u,
\end{multline}
and
\begin{multline}\label{dvweak}
\frac{\partial \mathcal{LF}_{\varepsilon}}{\partial \tilde{v}}(u_0,\tilde{v}_0,\phi_0)[\tilde{v}]=
-b\int_{\Omega}\psi'_{\varepsilon}(v_0)|\nabla u_0|^2\tilde{v}+\frac{c^2}{\varepsilon}\int_{\Omega}
(-V'(v_0)\tilde{v})
+\\2\varepsilon\int_{\Omega}\nabla\tilde{v}_0\cdot\nabla \tilde{v}
-\int_{\Omega}\psi'_{\varepsilon}(v_0)\tilde{v}\nabla u_0\cdot\nabla \phi_0,
\end{multline}
and, finally,
\begin{equation}\label{dpweak}
\frac{\partial \mathcal{LF}_{\varepsilon}}{\partial \phi}(u_0,\tilde{v}_0,\phi_0)[\phi]=
\int_{\Omega}\psi_{\varepsilon}(v_0)\nabla u_0\cdot\nabla \phi-\int_{\partial\Omega}f_{\varepsilon}\phi.
\end{equation}

Then the resulting optimality system is the following. We look for critical points, or better minimizers,
of $\mathcal{F}_{\varepsilon}$, or,
equivalently, of $\mathcal{MF}_{\varepsilon}(u,\tilde{v})$ subject to the constraint
$u=\mathcal{H}_{\varepsilon}(\tilde{v})$.
We use a gradient method, whose algorithm is divided into steps. A completely analogous method may be used for finding minimizers of
$\mathcal{G}_{\varepsilon}$.

\subsubsection{Step 0: initialization.}
We initialize the algorithm by putting $k=0$ and choosing an initial guess $\tilde{v}_0\in W(\Omega)$ such that
$0\leq \tilde{v}_0\leq 1$ almost everywhere. We observe that taking 
$\tilde{v}_0 \equiv 0$ (that is $v_0\equiv 1$) is not a good choice because this is a critical point of the functional
$\mathcal{F}_{\varepsilon}$, thus the gradient method fails in this case.

\subsubsection{Step 1: finding $\boldsymbol{u_k}$.}
We solve
\begin{equation}\label{dpstrong}
\left\{\begin{array}{ll}
\mathrm{div}(\psi_{\varepsilon}(v_k)\nabla u_k)=0 &\text{in }\Omega\\
\psi_{\varepsilon}(v_k)\nabla u_k\cdot\nu=f_{\varepsilon} &\text{on }\partial\Omega\\
\int_{\gamma}u_k=0, &
\end{array}\right.
\end{equation}
that is we look for $u_k\in W^{1,2}_{\gamma}(\Omega)$
such that
\begin{equation}\label{dpweak2}
\int_{\Omega}\psi_{\varepsilon}(v_k)\nabla u_k\cdot\nabla \phi-\int_{\partial\Omega}f_{\varepsilon}\phi=0\quad
\text{for any }\phi\in W^{1,2}(\Omega).
\end{equation}
We notice that $u_k=\mathcal{H}_{\varepsilon}(\tilde{v}_k)$ and
$u_k$ actually belongs to
$W_{\gamma}^{1,q(\varepsilon)}(\Omega)$.
By \eqref{lagrangian} and
by \eqref{dpweak}, we have that for any $\tilde{\phi}\in W^{1,2}(\Omega)$
$$\mathcal{LF}_{\varepsilon}(u_k,\tilde{v}_k,\tilde{\phi})=\mathcal{MF}_{\varepsilon}(u_k,\tilde{v}_k)=
\mathcal{F}_{\varepsilon}(\tilde{v}_k)\quad\text{and}\quad
\frac{\partial \mathcal{LF}_{\varepsilon}}{\partial \phi}(u_k,\tilde{v}_k,\tilde{\phi})=0.$$

\subsubsection{Step 2: finding $\boldsymbol{\phi_k}$.}
We solve
the following boundary value problem
\begin{equation}\label{dustrong}
\left\{\begin{array}{ll}
\mathrm{div}(\psi_{\varepsilon}(v_k)\nabla \phi_k)=-\mathrm{div}(2b\psi_{\varepsilon}(v_k)\nabla u_k) &\text{in }\Omega\\
\psi_{\varepsilon}(v_k)\nabla \phi_k
\cdot\nu=-\displaystyle{\frac{2a}{\varepsilon^{q_1}}(u_k-g_{\varepsilon})\chi_{\gamma}}
 &\text{on }\partial\Omega\\
\int_{\gamma}\phi_k=0. &
\end{array}\right.
\end{equation}
Here $\chi_{\gamma}$ denotes the characteristic function of $\gamma$, that is
$$(u_k-g_{\varepsilon})\chi_{\gamma}=\left\{\begin{array}{ll}
(u_k-g_{\varepsilon})&\text{on }\gamma\\
0&\text{on }\partial\Omega\backslash \gamma.
\end{array}\right.$$

The weak formulation of \eqref{dustrong} is looking for $\phi_k\in W^{1,2}_{\gamma}(\Omega)$ such that
\begin{multline}\label{duweak2}
\int_{\Omega}\psi_{\varepsilon}(v_k)\nabla \phi_k\cdot \nabla
u=\\-
2b\int_{\Omega}\psi_{\varepsilon}(v_k)\nabla u_k\cdot\nabla u-\frac{2a}{\varepsilon^{q_1}}
\int_{\gamma}(u_k-g_{\varepsilon})u\quad\text{for any }u\in W^{1,2}(\Omega).
\end{multline}
Such a solution $\phi_k$ exists and is unique.
Then $\mathcal{LF}_{\varepsilon}(u_k,\tilde{v}_k,\phi_k)=\mathcal{MF}_{\varepsilon}(u_k,\tilde{v}_k)=
\mathcal{F}_{\varepsilon}(\tilde{v}_k)$ and, by \eqref{duweak},
$$\frac{\partial \mathcal{LF}_{\varepsilon}}{\partial \phi}(u_k,\tilde{v}_k,\phi_k)=0\quad\text{and}\quad
\frac{\partial \mathcal{LF}_{\varepsilon}}{\partial u}(u_k,\tilde{v}_k,\phi_k)=0.$$

\subsubsection{Step 3: computing the gradient and updating $\boldsymbol{v_k}$.}
We compute the differential of $\mathcal{F}_{\varepsilon}$ at the point $\tilde{v}_k$.
We observe that if $u=\mathcal{H}_{\varepsilon}(\tilde{v})$, then for any $\tilde{\phi}\in W^{1,2}(\Omega)$ we have
$$\mathcal{F}_{\varepsilon}(\tilde{v})=\mathcal{MF}_{\varepsilon}(\mathcal{H}_{\varepsilon}(\tilde{v}),\tilde{v})=
\mathcal{LF}_{\varepsilon}(\mathcal{H}_{\varepsilon}(\tilde{v}),\tilde{v},\tilde{\phi}).$$
Therefore, since $u_k=\mathcal{H}_{\varepsilon}(\tilde{v}_k)$, and if we pick $\tilde{\phi}=\phi_k$, then
$$D\mathcal{F}_{\varepsilon}(\tilde{v}_k)=\frac{\partial \mathcal{LF}_{\varepsilon}}{\partial \tilde{v}}(u_k,\tilde{v}_k,\phi_k).$$
We conclude that, by \eqref{dvweak}, we have for any $\tilde{v}\in W_{p(\varepsilon)}(\Omega)$
\begin{multline}\label{dvweak0}
D\mathcal{F}_{\varepsilon}(\tilde{v}_k)[\tilde{v}]=
-b \int_{\Omega}\psi'_{\varepsilon}(v_k)|\nabla u_k|^2\tilde{v}+\frac{c^2}{\varepsilon}\int_{\Omega}
(-V'(v_k)\tilde{v})
+\\2\varepsilon\int_{\Omega}\nabla\tilde{v}_k\cdot\nabla \tilde{v}
-\int_{\Omega}\psi'_{\varepsilon}(v_k)\tilde{v}\nabla u_k\cdot\nabla \phi_k.
\end{multline}

Let us now consider the space $W_{p(\varepsilon)}(\Omega)$. We recall that
either $W_{p(\varepsilon)}(\Omega)=W_2(\Omega)$ (if $N=2$), that is $W_{p(\varepsilon)}(\Omega)$ is a Hilbert space,
or $W_{p(\varepsilon)}(\Omega)$ is a strictly convex real reflexive Banach space (if $N>2$). In either cases, if
$W=W_{p(\varepsilon)}(\Omega)$,
we fix an operator $T:W^{\ast}\to W$ such that for any $w^{\ast}\in W^{\ast}$, we have
$$\langle w^{\ast},T(w^{\ast})\rangle= \|w^{\ast}\|^2\quad\text{and}\quad \|T(w^{\ast})\|=\|w^{\ast}\|,$$
where $\langle \cdot,\cdot\rangle$ is the usual duality between $W^{\ast}$ and $W$.
We may choose $T$ as the duality mapping from $W^{\ast}$ into $W^{\ast\ast}=W$.
If $W$ is a Hilbert space and we also identify $W^{\ast}$ with
$W$, then $T$ is actually the identity.
See, for instance, \cite[Section~42.6]{Zei}. Let us call $T_{\varepsilon}$ the corresponding operator for $W_{p(\varepsilon)}(\Omega)$.

For a positive constant $t_k$, we then update $\tilde{v}_k$ by setting
$$\hat{v}_{k+1}=\tilde{v}_k-t_kT_{\varepsilon}(D\mathcal{F}_{\varepsilon}(\tilde{v}_k)).$$

We observe the following. If $D\mathcal{F}_{\varepsilon}(\tilde{v}_k)=0$,
then $(u_k,\tilde{v}_k,\phi_k)$ is a critical point of $\mathcal{LF}_{\varepsilon}$ and $\tilde{v}_k$ is a critical point of $\mathcal{F}_{\varepsilon}$ and the algorithm comes to a stop.
Otherwise, provided $t_k$ is small enough, an easy computation shows that
$\mathcal{F}_{\varepsilon}(\hat{v}_{k+1})<\mathcal{F}_{\varepsilon}(\tilde{v}_k)$.

\subsubsection{Step 4: normalization and finding $\boldsymbol{\tilde{v}_{k+1}}$.}
We normalize $\hat{v}_{k+1}$ by truncation as follows. We set $\tilde{v}_{k+1}=(\hat{v}_{k+1}\wedge 1)\vee 0$. In such a way we obtain that $\tilde{v}_{k+1}\in W_{p(\varepsilon)}(\Omega)$ and
$0\leq \tilde{v}_{k+1}\leq 1$ almost everywhere in $\Omega$.

Let us note that, by our hypotheses, such a truncation does not increase the value of the functional,
in fact for any $\hat{v}\in W_{p(\varepsilon)}(\Omega)$, if 
$\tilde{v}=(\hat{v}\wedge 1)\vee 0$, then
$$\mathcal{F}_{\varepsilon}(\tilde{v})\leq\mathcal{F}_{\varepsilon}(\hat{v}).$$
Therefore, we have found that either $D\mathcal{F}_{\varepsilon}(\tilde{v}_k)=0$,
and the algorithm stops, or, otherwise,
provided $t_k$ is small enough,
$\mathcal{F}_{\varepsilon}(\tilde{v}_{k+1})<\mathcal{F}_{\varepsilon}(\tilde{v}_k)$.

Once we have computed $\tilde{v}_{k+1}$, we iterate the algorithm by going back to Step~1.

\section{Numerical experiments}\label{numsec}

The data for the numerical experiments are generated by solving Laplace equation numerically on an domain with certain prescribed defects (cracks or cavities).
We solve the Neumann problem with given flux on the boundary of the computational domain, and read off the corresponding Dirichlet data to 
get a feasible pair of Neumann and Dirichlet boundary data on a discrete set of measurement points on the boundary from which the defect has to be reconstructed. As input fluxes we choose pairs of plus-shaped current profiles with opposite sign located at two different sides of the rectangular computational domain. The Laplace equation is solved on a very fine irregular grid using linear finite elements. The boundary data are genuinely defined on the unevenly distributed nodal points of elements on the boundary and are interpolated onto a much courser regular grid of measurement points. When experimenting with noisy input data, both boundary values are contaminated by adding Gaussian distributed artificial noise to the data, usually with different noise levels for $f = \frac{\partial u}{\partial \nu}|_\gamma$ and $g = u|_\gamma$.

For the numerical implementation of step 1 in the algorithm described in the previous section (that is the numerical solution of equation \eqref{dpweak2} for $u_k$ with given $v_k$ and prescribed $f_\varepsilon$), we also use linear finite elements for the discretization of $u_k$. In contrast to the data generation routine, we discretize the potential on a \emph{regular, structured grid} which is usually much coarser than the grid used for the data generation. Later on, we shall assume that the phase-field $v_k$ is also an element in the space of piecewise linear functions on the same underlying regular grid as for $u_k$. For the assembling of the stiffness matrix for  \eqref{dpweak2}, however, we replace the phase-field $v_k$ by its $L^2$-projection onto the space of functions which are piecewise constant on the triangles of the finite element space. A completely analogous procedure is applied for the solution of the adjoint equation \eqref{duweak2} described in step 2 for the adjoint variable $\phi_k$. Note that both systems share the same stiffness matrix and that the right-hand side of \eqref{duweak2} can be easily assembled using a slightly modified stiffness matrix. We shall use up to six different Cauchy data-sets for the reconstruction of the defect. The data-sets correspond to all possible combinations of pairs of electrodes where each electrode is located on a different side of the computational rectangle. We can use the same factorization of the stiffness matrix for all different right-hand sides of \eqref{dpweak2} and \eqref{duweak2}.

The calculation of the descent direction for the cost functional as described in step 3 requires another solution of an elliptic boundary value problem for the variable $\delta \tilde{v}_k = T_{\varepsilon}(D\mathcal{F}_{\varepsilon}(\tilde{v}_k))$. As mentioned above, the update $\delta \tilde{v}_k$ is discretized using linear triangular elements on a regular grid. To find $\delta \tilde{v}_k$ we have to solve an elliptic equations with system matrix defined by a discretization of the operator $T:W^\ast \to W$. In our 2-dimensional test examples, we always set $W = W^{1,2}(\Omega)$ and
for any $w^{\ast}\in W^{\ast}$ we set $T(w^{\ast})=v$ where $v$ solves in a weak sense $v-c\Delta v=w^{\ast}$
with some parameter $c>0$ and homogeneous Dirichlet boundary conditions. The choice of Dirichlet boundary conditions is motivated by the desire to keep the phase-field constantly at the value 1 on the boundary. The assembling of the right-hand side of the equation for $\delta \tilde{v}_k$ is done by evaluating \eqref{dvweak0} for piecewise linear in all bases functions $\tilde{v}$. 

The projection required in step 4 is easily implemented for piecewise linear functions by thresholding the nodal values. Moreover, a suitable step-length for the update of the phase-field is found using an Armijo-type line search. We use a maximum number of five reduction steps for the correction of the step-length. Since each evaluation of the cost functional requires one solution of the state equation, we try to steer the step-size modification in a rather conservative way. 

\begin{figure}[htb]
\centering
\includegraphics[width = \textwidth/2]{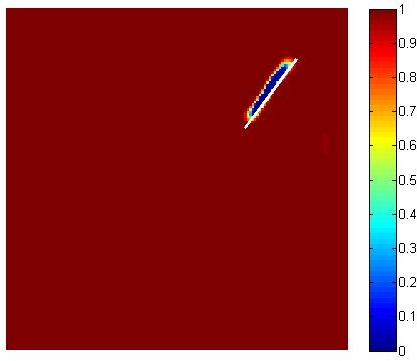}
\caption{Reconstruction of a small linear crack with noise-free data.}
\label{fig1}
\end{figure}

Within this setup, the following numerical experiments have been performed. For all experiments, the phase-field parameter $\varepsilon$ was decreased in several steps from an initial value of $\varepsilon = 2 \cdot 10^{-4}$ down to $\varepsilon = 1 \cdot 10^{-6}$ for the single-well potential and to $\varepsilon = 2 \cdot 10^{-6}$ for the double-well case. We run 2500 iterations of our algorithm in the single-well case and 1000 in the double-well case.
Figure \ref{fig1} shows the final phase-field together with the linear crack (as a white line) which was used for the data generation. We use all six available data-sets with electrode positions on (up/down), (left/right), (down/left), (up/left), (down/right), and (up/right) sides of the rectangle for the reconstruction and set the noise-level to zero. In this simple situation where the crack is located rather close to the boundary we obtain very good reconstruction of the crack location with the single-well approximation.

\begin{figure}[htb]
\centering
\mbox{\includegraphics[width = .45\textwidth]{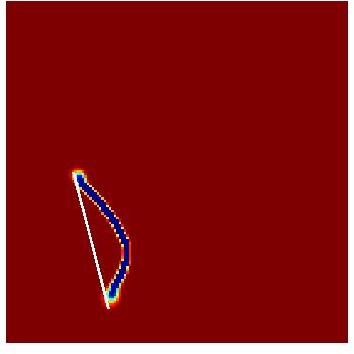}
\includegraphics[width = .45\textwidth]{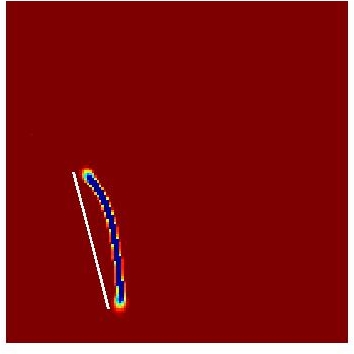}}
\caption{Comparison of reconstructions from 3 and 6 measurements.}
\label{fig2}
\end{figure}

In Figure \ref{fig2} it is shown a comparison between reconstructions using 3 measurements (left image) with electrode positions on (left/right), (left/up), and (right/up) edges and 6 measurements (right figure), again in the single-well case. It is notable that in the reconstruction with 3 data-sets the crack tips are accurately identified but the reconstructed crack is strongly curved which is probably due to the fact that we have no electrode located on the lower edge of the computational domain. In contrast the overall geometrical shape of the crack is reconstructed much better with 6 data-sets but the position of the crack tips is less accurate. In these two simulations we added one percent of normally distributed noise to Neumann and Dirichlet data. 

\begin{figure}[htb]
\centering
\mbox{\includegraphics[width = .45\textwidth]{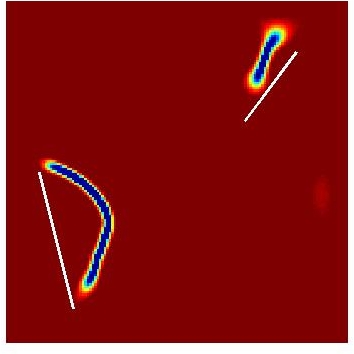}
\includegraphics[width = .45\textwidth]{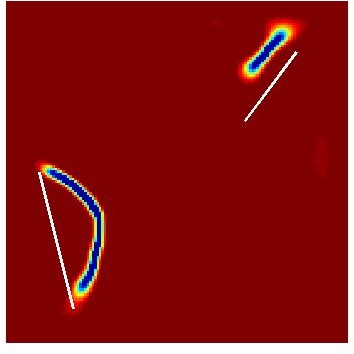}}
\caption{Comparison with different noise levels.}
\label{fig3}
\end{figure}

Figure \ref{fig3} shows results for a situation with two cracks and different noise levels. Here we fixed the noise-level for the Neumann data to 1\% for both experiments whereas the Dirichlet data were contaminated with 1\% (left image) and 5\% (right image) of noise. We used three measurements (left/right), (left/up), (right/up) and the single-well potential. There is no big difference in the quality of the reconstructions. In both cases the placement of the smaller crack in the upper right corner is inaccurate and the larger crack in the lower left corner is curved. Nonetheless the convergence of the algorithm is not heavily effected by the presence of (moderately strong) noise and the reconstructions are stable. 

\begin{figure}[htb]
\centering
\mbox{\includegraphics[width = .32\textwidth]{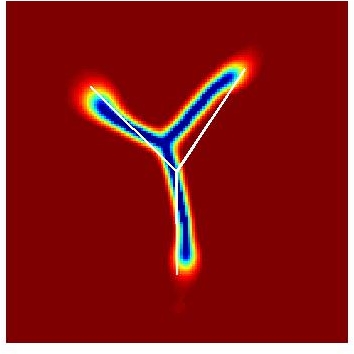}
\includegraphics[width = .32\textwidth]{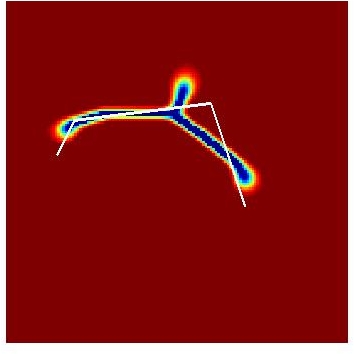}
\includegraphics[width = .32\textwidth]{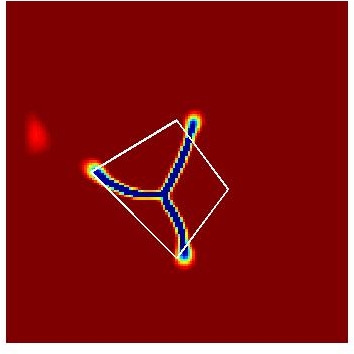}}
\caption{Dendrite-like reconstructions with single-well potential.}
\label{fig4}
\end{figure}

The next series of experiments presented in Figure \ref{fig4} shows the tendency of the single-well based algorithm to produce dendrite-like structures. In fact,
the dendrite-shaped crack in the leftmost image is reconstructed quite well. The polygonal crack in the middle image is approximated by a cloth-hanger like structure which has a satisfactory data fit with a shorter overall length than the polygonal curve. Finally the cavity in the rightmost image is approximated by a one-dimensional structure which looks roughly like the skeleton of the cavity. In all these three experiments noise level is 1\% for Neumann data and 5\% for Dirichlet data and the three measurements (left/right), (left/up), (right/up) are used.

\begin{figure}[htb]
\centering
\mbox{\includegraphics[width = 0.45\textwidth]{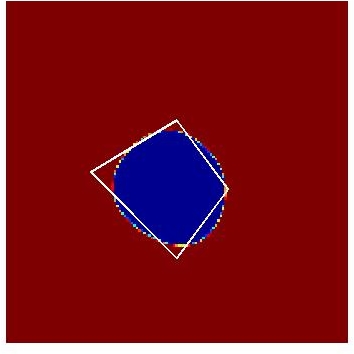}
\includegraphics[width = 0.45\textwidth]{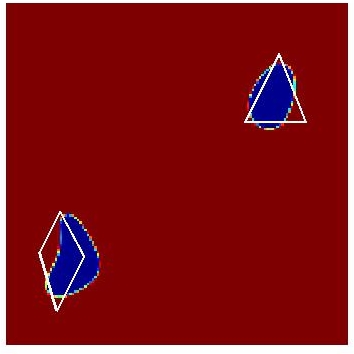}}
\caption{Reconstructions of cavities with double-well potential.}
\label{fig5}
\end{figure}

Figure \ref{fig5} shows reconstructions obtained by using the double-well approximation. As expected, the phase-field approximates the characteristic functions of one cavity (left image) and two cavities (right image).
In these two tests noise level is 1\% for Neumann data. In the left image noise level for Dirichlet data is 5\% and 
the three measurements (left/right), (left/up), (right/up) are used. In the right image noise level for Dirichlet data is 1\% but only one measurement, namely (left/right), is used.
The overall location of the cavities is satisfactory, but the lower left quadrilateral is approximated by a non-convex shape. In this respect the experiment with the double-well potential for two cavities resembles the results shown in Figure \ref{fig2} where the lower left crack also has a strong tendency  to bend inward.

For our final numerical experiment, documented in Figure \ref{fig5}, the double-well approach was used for the reconstruction of one-dimensional defects like the polygonal crack shown in the left image and the star-shaped crack shown on the right-hand side of the figure. In both cases the defect is approximated by a two dimensional structure. An interesting feature is the occurrence of a self-intersection of the boundary curve of the reconstructed defect in the case of the star-shaped crack.
Also in these two final tests, noise level is 1\% for Neumann data and 5\% for Dirichlet data and 
the three measurements (left/right), (left/up), (right/up) are used.
\begin{figure}[htb]
\centering
\includegraphics[width = 0.45\textwidth]{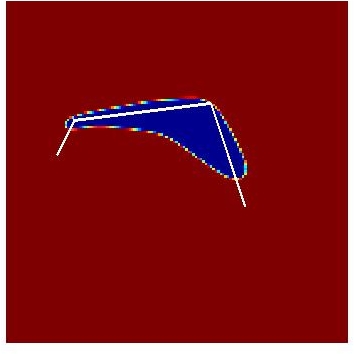}
\includegraphics[width = 0.45\textwidth]{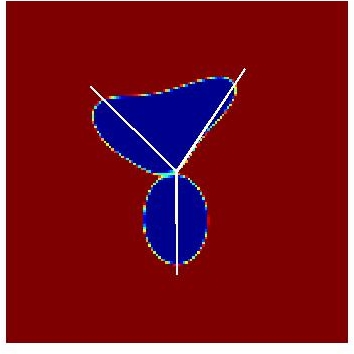}
\caption{Reconstructions of cracks with double-well potential.}
\label{fig6}
\end{figure}

As a conclusion we can state that both algorithms give reconstructions of the defects with a satisfactory accuracy for an exponentially ill-posed problem. The algorithms show a quite stable behaviour in the presence of data noise. The single-well and double-well models develop the types of structures for which they are designed (one-dimensional for the single-well and two dimensional for the double-well potential), so the single-well approach approximates cavities by dendrites and the double-well approach approximates cracks by cavities. The double-well approach looks more stable with respect to noise, is slightly less sensitive with respect to the adjustment of the phase-field parameter $\varepsilon$ and usually needs less iterations for convergence. This may be in accordance with the theory, in fact for the double-well case a convergence analysis is proved, whereas the single-well model we use is a modification of the one for which we have convergence results. Finally, it turned out to be important to update the phase-field parameter $\varepsilon$ adaptively during the algorithm. If the parameter $\varepsilon$ is chosen too small initially or decreased too fast, sharp interfaces develop too early, sometimes at incorrect locations, and the algorithm is not able to move well established interfaces to other locations. On the other hand, if the parameter $\varepsilon$ is decreased too much, the term containing the potential might prevail and not well established defects, usually the smaller ones, may disappear.

\end{document}